\documentclass[reqno]{amsart}
\usepackage{amsmath,amsfonts,amssymb,graphics,graphicx,latexsym,amscd,subfigure,hyperref}
\usepackage[usenames]{xcolor}
\usepackage[retainorgcmds]{IEEEtrantools}

\DeclareMathAlphabet{\mathantt}{OT1}{antt}{li}{it}
\DeclareMathAlphabet{\mathpzc}{OT1}{pzc}{m}{it}

\newtheorem{thm}{Theorem}[section]

\newtheorem{question}[thm]{Question}
\newtheorem{lemma}[thm]{Lemma}
\newtheorem{defn}[thm]{Definition}

\newcommand{\bbR}{\mathbb{R}}
\newcommand{\bbT}{\mathbb{T}}

\newcommand{\bbZ}{\mathbb{Z}}
\newcommand{\ip}{\cdot}
\newcommand{\cp}{P}
\newcommand{\s}{g}
\newcommand{\im}{\sqrt{-1}}

\newcommand{\hess}{\operatorname{Hess}}

\newcommand{\Lap}{\triangle}

\newcommand{\tr}{\operatorname{Tr}}

\newcommand{\action}{\psi}

\def\det{\mathop{\rm det}\nolimits}

\def\sgn{\mathop{\rm sgn}\nolimits}

\def\tr{\mathop{\rm tr}\nolimits}

\begin{document}

\title[]{Recovering $S^1$-invariant metrics on $S^2$ from the equivariant spectrum}

\author{Emily B. Dryden}
\address{Department of Mathematics, Bucknell University, Lewisburg, PA 17837, USA}
\email{ed012@bucknell.edu}
\author{Diana Macedo}
\address{Departamento de Matem\'{a}tica, Instituto Superior T\'{e}cnico, Av. Rovisco Pais, 1049-001 Lisboa, Portugal}
\email{dnmcd6@hotmail.com}
\author{Rosa Sena-Dias}
\address{Centro de An\'alise Matem\'atica, Geometria e Sistemas Din\^amicos, Departamento de Matem\'{a}tica, Instituto Superior T\'{e}cnico, Av. Rovisco Pais, 1049-001 Lisboa, Portugal}
\email{rsenadias@math.ist.utl.pt}
\thanks{EBD was partially supported by a grant from the Simons Foundation (210445 to Emily B. Dryden).  DM was partially supported by the Gulbenkian Foundation through a ``Novos talentos em Matem\' atica" grant. RSD was partially supported by FCT/Portugal through projects PEst-OE/EEI/LAOO9/2013, EXCL/MAT-GEO/0222/2012 and PTDC/MAT/117762/2010}

\date{}

\begin{abstract}
We prove an inverse spectral result for $S^1$-invariant metrics on $S^2$ based on the so-called asymptotic equivariant spectrum. This is roughly the spectrum together with large weights of the $S^1$ action on the eigenspaces. Our result generalizes an inverse spectral result from \cite{dgs3} concerning $S^1$-invariant metrics on $S^2$ which are invariant under the antipodal map. We use higher order terms in the asymptotic expansion of a natural spectral measure associated with the Laplacian and the $S^1$ action.
\end{abstract}

\maketitle


\section{Introduction}\label{sec:intro}

Does the spectrum of the Laplacian on {a} compact Riemannian manifold determine the Riemannian manifold? The answer is known to be no.  John Milnor \cite{m} constructed the first such example, producing two metrics on flat tori of dimension 16 which have the same spectrum, but are not globally isometric. There are now a plethora of examples of non-isometric isospectral manifolds, constructed using sophisticated methods and tools.  On the other hand, there are also many \emph{spectral invariants}, quantities associated with the Riemannian manifold  {that} are determined by the spectrum. Dimension and volume are perhaps the most famous of these spectral invariants, but many more are known. Given these extremes of examples and spectral invariants, it is natural to ask if there is a special class of Riemannian manifolds for which the spectrum determines the Riemannian structure. In this spirit, Miguel Abreu formulated the following question.
\begin{question}\label{q:Abreu}
Let $X$ be a toric manifold endowed with a toric K\"ahler structure. Does the spectrum of the Laplacian on the resulting Riemannian manifold determine the underlying manifold?
\end{question}

Toric manifold{s} are rather manageable objects.  They are classified by convex polytopes of so-called Delzant type and the question above can be reformulated in terms of recovering the Delzant polytope of a toric manifold from the spectrum of a toric K\"ahler metric on it. (See \cite{dgs1} for some basic background on toric manifolds and Delzant polytopes; see \cite{g} for a more complete account of the K\"ahler geometry of toric manifolds.)
In the setting of toric orbifolds, a partial answer to Question \ref{q:Abreu} using the \emph{asymptotic equivariant spectrum} is known; the asymptotic equivariant spectrum adds information about the large weights of the  {$\bbT^n$} action on the eigenspaces to the usual spectrum.
\begin{thm}\cite{dgs3}
The asymptotic equivariant spectrum of a {\em generic} toric orbifold endowed with a toric K\"ahler metric determines the toric orbifold up to equivariant biholomorphism.
\end{thm}

Since the asymptotic equivariant spectrum contains more information than the spectrum, and toric orbifolds are more ``flexible'' than toric manifolds, the above theorem does not answer Abreu's question. However, it is also possible to obtain spectral rigidity for the {\em metric}, at least in a simple case.  Note that a toric K\"{a}hler metric can be entirely described in terms of a function $\s$ known as its symplectic potential.  {We will consider $S^2$, suitably normalized} so that its symplectic potential is a function $\s$ defined on $(-1,1)$. Given the symplectic potential $\s$, the metric can be written as
\begin{equation}\label{metric_in_action_angle}
\ddot{\s}dx \otimes dx + \frac{d \theta \otimes d \theta}{\ddot{\s}}.
\end{equation}

\begin{thm}\cite{dgs3}\label{thm:DGS3metric}
The asymptotic equivariant spectrum of an $S^1$-invariant metric on $S^2$ with symplectic potential $\s$ determines the metric if $\ddot{\s}$ is even and convex.
\end{thm}

Our goal is to extend Theorem \ref{thm:DGS3metric} to more general $S^1$-invariant metrics on $S^2$.  {In particular, we will consider metrics for which $\ddot{\s}$ is a single well.}  A single well on $(-1,1)$ will mean a function $v:(-1,1) \rightarrow \bbR$ which is decreasing on $(-1,0)$ and increasing on $(0,1)$. The unique minimum at $0$ is the well.  We will prove the following theorem.
\begin{thm}\label{main}
The asymptotic equivariant spectrum of an $S^1$-invariant metric on $S^2$ with symplectic potential $\s$ determines the metric if $\ddot{\s}$ is a single well.
\end{thm}
 {In other words,} we show that the asymptotic equivariant spectrum determines an $S^1$-invariant metric on $S^2$ within the class of metrics  {for which $\ddot{\s}$} is a single well. It would be interesting to know if the asymptotic equivariant spectrum determines if an $S^1$-invariant metric on $S^2$ has such a  {symplectic} potential. 

Although this result concerns only $S^2$, results on spectral rigidity for metrics are known to be difficult to obtain. One of the reasons we are able to prove such a result in this setting is the parametrization of $S^1$-invariant metrics on $S^2$ via the symplectic potential.  This framework was also used by Miguel Abreu and Pedro Freitas \cite{af} to investigate upper bounds for the invariant eigenvalues of the Laplace operator defined by these $S^1$-invariant metrics on $S^2$. 

{Related} inverse spectral result{s} for $S^1$-invariant metrics on $S^2$ w{ere} obtained by {Jochen Br{\"u}ning and Ernst Heintze \cite{bh} and by} Steve Zelditch {\cite{z}}. {Zelditch's} result concerns $S^1$-invariant metrics on $S^2$ that can be obtained as metrics on surfaces of revolution in $\mathbb{R}^3$ with the euclidean metric. We refer to such metrics as metrics of revolution. The main theorem in \cite{z} says that generic,  {``simple,''} analytic metrics of revolution are determined by their spectrum.  Not all $S^1$-invariant metrics on $S^2$ are metrics of revolution on $\bbR^3$ (see \cite{kw} for obstructions  {and \cite[\S4]{af} for discussion}), but for an $S^1$-invariant metric on $S^2$  {that} is a metric of revolution the  {``simple''} condition is {implied by} our single well condition. To be more precise, a surface of revolution in $\bbR^3$ can be obtained by rotating a curve $(0,p(t),q(t)))$  {around} the $z$-axis where we assume that $t$ is the arc-length parameter. {The metric of revolution described above is given by
$dt^2+p(t)^2d\theta^2$, 
where $\theta$ is the angle of revolution. By comparing with (\ref{metric_in_action_angle}) we see that 
$$
t(x)=\int_{-1}^x\sqrt{\ddot{\s}(s)}ds\qquad \text{\ \ and \ \ } \qquad
p(t)^2=\frac{1}{\ddot{\s}(x)}.
$$
Therefore local minima of $\ddot{\s}$ (at $x$) correspond to local maxima of $p$ (at $t(x)$). The ``simple" condition corresponds to $p$ having a single local maximum and the claim follows.} The main theorem in \cite{z} is proved using wave trace techniques. {Br{\"u}ning and Heintze considered metrics of revolution that are also symmetric with respect to reflection in the $xy$-plane.  They showed that mirror symmetric surfaces of revolution are determined by their spectra, or equivalently, by their $S^1$-invariant spectra.}

The techniques we use to prove  {Theorem \ref{main}} are inspired by the techniques of Victor Guillemin and Zuoqin Wang in \cite{gw}. Namely, we give an explicit inductive formula giving higher order semi-classical spectral invariants for toric manifolds based on the asymptotic expansion of a spectral measure associated to the toric metric. This spectral measure was introduced in \cite{dgs3} but there only the highest term in its asymptotic expansion was derived and used. The formula for higher order invariants turns out to be surprisingly explicit in the context of toric manifolds and is interesting in its own right. In the case of $S^2$ it becomes even simpler and we are able to derive some new spectral invariants which we use to prove Theorem \ref{main}.

One can't help but wonder if semi-classical spectral invariants of even higher order on $S^2$ would help recover symplectic potentials with a finite number of wells from the asymptotic equivariant spectrum. It would also be interesting to study these invariants in higher dimensional toric manifolds.

The paper is organized as follows: in \S\ref{asymp_expansion} we derive a general formula for the asymptotic expansion of a spectral measure associated to the Laplacian of a toric K\"ahler metric on a toric manifold. The formula simplifies in the case of $S^2$ and we give two non-trivial terms in the asymptotic expansion. The first term was given and used in \cite{dgs3}. Section \ref{proof} is devoted to the proof of Theorem \ref{main}. We use both spectral invariants derived in \S\ref{asymp_expansion}. \\

\noindent \textbf{Acknowledgements.}  We would like to thank Victor Guillemin for his encouragement and interest in this project.  Emily Dryden appreciates the hospitality of the Mathematics Department at Instituto Superior T\'{e}cnico during the preparation of this paper. Diana Macedo and Rosa Sena-Dias would like to thank the Gulbenkian Foundation and in particular the ``Novos talentos em Matem\' atica" program organizers for support and for creating a great environment for mathematical research. We  {are grateful to} the referees for carefully reading our manuscript and for  {their} helpful suggestions.


\section{Asymptotic expansion of the spectral measure}\label{asymp_expansion}

Let $X^{2n}$ be a toric manifold with moment polytope $\cp$. Suppose $X$ is endowed with a toric K\"ahler metric defined by a symplectic potential $\s: \operatorname{Int}(\cp)\rightarrow \bbR$, i.e., a metric of the form 
\[
ds^2=\ddot{\s}dx \otimes dx + \frac{d \theta \otimes d \theta}{\ddot{\s}},
\]
(see \cite{a} for the definition of symplectic potential and the associated metric). The torus $\bbT^n$ acts on $X$ by isometries, and we let $\action:\bbT^n \rightarrow \text{Iso}(X,ds^2)$ denote the action. The metric on $X$ has a Laplace operator associated to it. We are interested in the question: Does the spectrum of this Laplacian determine the metric, or equivalently, the symplectic potential $\s$? In fact, the torus action induces a representation on each eigenspace  {that} splits according to weights in the Lie algebra of the torus, suitably identified with $\bbR^n$. We will sometimes write $\action$ for the map from {$\bbR^n/(2\pi\bbZ)^n$} to $\text{Iso}(X,ds^2)$ where we identify $\bbT^n$ with $\bbR^n/{(2\pi\bbZ)}^n$ via the exponential map. The {\em equivariant spectrum} of  {($X, ds^2$) is the set of eigenvalues of the Laplacian and}, for each eigenvalue, the list of weights of the torus representation on the corresponding eigenspace. The eigenvalues and weights are counted with multiplicities. We say that a toric K\"ahler metric or any quantity associated with it is \emph{espectrally determined} if it is determined by its equivariant spectrum (see \cite{dgs1} for more details about the equivariant spectrum).
\begin{defn}
Let $K>0$. The $K$-equivariant spectrum of a toric K\"ahler manifold is the set of eigenvalues of its metric Laplacian and, for each eigenvalue, the list of weights of the torus representation on the corresponding eigenspace whose norms are larger than $K$.  The eigenvalues and weights are counted with multiplicities.
\end{defn}

We will say that a toric K\"ahler metric or any quantity associated with it is aespectrally determined if there exists a $K>0$, possibly very large, for which it is determined by the $K$-equivariant spectrum. 

Let $\alpha$ be a generic element in $\bbR^n$, where ``generic'' means that $\alpha$ is a regular value of the moment map $\Phi: T^*X \rightarrow \bbR^n$ arising from the lift of $\action$ to a Hamiltonian action on $T^*X$.  {Moreover, let $\hbar \in \bbR^+$ be such that $\frac{1}{\hbar}\in \bbZ$.  We wish to define a spectral measure that ``counts" $\frac{\alpha}{\hbar}$-equivariant eigenfunctions in a certain sense.  Let $\rho\in \mathcal{C}_0^\infty(\bbR)$ and consider the spectral measure $\mu_{\frac{\alpha}{\hbar}}$ defined by
\begin{equation}\label{eqn:measure}
\mu_{\frac{\alpha}{\hbar}}(\rho)=\int_{\bbT^n}e^{-\frac{{\im}\theta \cdot\alpha}{\hbar}}\tr(\action(\theta)_*\rho(\hbar^2\Lap)) \mathrm{d}\theta.
\end{equation}
Note that since $\frac{1}{\hbar}\in \bbZ$, $e^{-\frac{{\im}\theta \cdot\alpha}{\hbar}}$ is well defined.  Also, recall that given any pseudo-differential semi-classical operator $\mathcal{P}_\hbar$ on a manifold $X$ and a diffeomorphism $F$ of $X$, $\tr(F_* \mathcal{P_\hbar})$ is defined using the Schwartz kernel $\mathcal{K}$ of $\mathcal{P_\hbar}$ via the formula
$$
\tr (F_*\mathcal{P_\hbar})=\int_X \mathcal{K}(F(x),x)dx,
$$
when such an integral converges. It turns out that $\mu_{\frac{\alpha}{\hbar}}$ may also be written as 
$$
\mu_{\frac{\alpha}{\hbar}}(\rho)=\sum_{\lambda \in E_{\alpha,\hbar}}\rho(\lambda)n(\lambda, \alpha),
$$
where $E_{\alpha,\hbar}$ is the set of eigenvalues of $\hbar^2\Lap$ admitting an equivariant eigenfunction of weight $\frac{\alpha}{\hbar}$ and $n(\lambda,\alpha)$ is the dimension of the subspace of $\frac{\alpha}{\hbar}$-equivariant functions of the $\lambda$-eigenspace of $\hbar^2\Lap$; this formulation of $\mu_{\frac{\alpha}{\hbar}}$ shows why we may interpret our spectral measure as ``counting'' equivariant eigenfunctions of weight $\frac{\alpha}{\hbar}$. See \cite[\S4]{dgs3} for more details.} 

As shown in \cite{dgs3}, the measure $\mu_{\frac{\alpha}{\hbar}}$ has an asymptotic expansion in powers of $\hbar$ as $\hbar$ tends to zero, and the terms in that expansion are aespectrally determined; the first term in the expansion is calculated explicitly.  We will give an algorithm to calculate terms in the expansion for toric manifolds. Our work is very much in the spirit of \cite{gw}, in which the authors describe an analogous algorithm for semi-classical operators in $\bbR^n$. Because toric manifolds admit an open dense set with global coordinates, we will use essentially the same techniques as in \cite{gw}.

In particular, we will use this algorithm to calculate the second nontrivial term in the asymptotic expansion for $S^2$ with an $S^1$-invariant metric. The first term of the expansion was used in \cite{dgs3} to show that, for $S^1$-invariant metrics on $S^2$, when the symplectic potential is convex and even it is aespectrally determined.  Using the first and second terms, we are able to extend the result to single well symplectic potentials. This is analogous to the situation in \cite{gw}, in which the authors extended a result of Yves Colin de Verdi\`ere \cite{cdv} on Schr\"{o}dinger operators to more general potentials.

We begin with the asymptotic expansion of $\mu_{\frac{\alpha}{\hbar}}$.
\begin{thm}\label{asymp_toric}
Let $X^{2n}$ be a toric manifold with a $\bbT^n$ action $\action$. Let $X$ be endowed with a toric K\"ahler metric whose symplectic potential is $\s$. Let $\alpha\in \bbR^n$ be generic. Then the measure $\mu_{\frac{\alpha}{\hbar}}$ defined in \eqref{eqn:measure} has an asymptotic expansion in powers of $\hbar$ given by
$$
\hbar^{-n}\sum_{k\geq0}\hbar^k\sum_{l\leq 2k}\int_{\cp\times \bbR^n}b_{k,l}^\alpha(u,\hat{u})\frac{1}{(\im)^l}\frac{d^l\rho}{d\sigma^l}\left( \hat{u}^t\hess^{-1}(\s)\hat{u}+\alpha^t\hess(\s)\alpha\right)\mathrm{d}u\mathrm{d}\hat{u};
$$
the functions $b^\alpha_{k,l}$ are given by $b^\alpha_{k}(u,\hat{u},t)=\sum_{l\leq2k}b^\alpha_{k,l}(u,\hat{u})t^l$ where $b^\alpha_{k}(u,\hat{u},t)$ are given recursively by $b^\alpha_0(u,\hat{u},t)=1$ and
\begin{IEEEeqnarray*}{l}
\frac{1}{\im}\frac{\partial b^\alpha_{1}}{\partial t}=2t\sum_{i,j=1}^n \s^{ij}\hat{u}_j \left( \hat{u}^t\frac{\partial\hess^{-1}(\s)}{\partial u_i}\hat{u}+\alpha^t\frac{\partial\hess(\s)}{\partial u_i}\alpha\right),\nonumber\\
\end{IEEEeqnarray*}
\begin{IEEEeqnarray*}{lCr}
\frac{1}{\im}\frac{\partial b^\alpha_{k}}{\partial t}&=&\frac{2}{\im}\sum_{i,j=1}^n\s^{ij}\hat{u}_j  \left(\frac{\partial}{\partial u_i}+\im t \left( \hat{u}^t\frac{\partial\hess^{-1}(\s)}{\partial u_i}\hat{u}+\alpha^t\frac{\partial\hess(\s)}{\partial u_i}\alpha\right)\right)b^{\alpha}_{k-1}\nonumber \\
&&-\sum_{i,j=1}^n \s^{ij}\left(\frac{\partial}{\partial u_i}+\im t \left( \hat{u}^t\frac{\partial\hess^{-1}(\s)}{\partial u_i}\hat{u}+\alpha^t\frac{\partial\hess(\s)}{\partial u_i}\alpha\right)\right)\cdot \nonumber \\
&&\left(\frac{\partial}{\partial u_j}+\im t \left( \hat{u}^t\frac{\partial\hess^{-1}(\s)}{\partial u_j}\hat{u}+\alpha^t\frac{\partial\hess(\s)}{\partial u_j}\alpha\right)\right) b^{\alpha}_{k-2}, 
\end{IEEEeqnarray*}
for all $k>1$, and $b_{k}^\alpha(u,\hat{u},0)=0$ for $k \geq 1$. Here $\s_{ij}$ denote the entries of $\hess(\s)$ and $\s^{ij}$ denote the entries of its inverse $\hess^{-1}(\s)$.
\end{thm}

The proof of this theorem is very similar to the proof of Theorem 5.1 in \cite{dgs3}. Since we are treating the case of the Laplace operator for a toric K\"{a}hler metric, rather than the more general case of a Riemannian metric on a manifold that admits an isometric action of some torus as in \cite{dgs3}, we have global coordinates on  {an} open dense set and some simplifications occur. For the convenience of the reader we will essentially give a complete proof. We need two main ingredients for this proof.

\begin{lemma}[Schwartz kernel asymptotic expansion]\label{lemma:schwartz}
With the setup and notation given above, $\rho(\hbar^2\Lap)$ is a semi-classical operator on $X$ with Schwartz kernel $K_{\rho,\hbar}$. In local coordinates $K_{\rho,\hbar}$ admits an asymptotic expansion in powers of $\hbar$:
$$
K_{\rho,\hbar}({\bf x},{\bf y})=(2\pi\hbar)^{-2n}\sum_{k\geq 0} {\hbar}^k\sum_{l\leq 2k}\int_{\bbR^{2n}}b_{k,l}({\bf y},{\bf \xi})\frac{1}{(\im)^l}\frac{d^l\rho}{d\sigma^l}(|{\bf \xi}|_\s^2({\bf x}))e^{\frac{\im({\bf x}-{\bf y})\ip {\bf \xi}}{\hbar}}\mathrm{d}{\bf \xi},
$$
where $|\cdot|_\s^2({\bf x})$ denotes the norm on the cotangent space over ${\bf x}$ given by the metric associated to $\s$ and where $b_{k,l}$ are given as follows. Let $b_{k}$ be defined recursively by $b_{0}({\bf x},{\bf \xi},t)=1$, $b_{k}({\bf x},{\bf \xi},0)=0$ for $k \geq 1$, and
\begin{equation}\label{inductionb_m_general}
\frac{1}{\im}\frac{\partial b_k}{\partial t}=\sum_{\substack{a=(a_1,\ldots, a_{2n})\in (\bbZ_+)^{2n} \\  |a|\geq 1}}\sum_{j+|a|=k}D^a_{{\bf \xi}}(|{\bf \xi}|_\s^2)Q_ab_j
\end{equation}
where
$$
D^a_\xi=\frac{1}{(\sqrt{-1})^{|a|}}\frac{\partial^{|a|}}{\partial \xi^a} 
 {\text{\ \  and \ \ }}
Q_a=\frac{1}{a!}\left(\frac{\partial}{\partial \bf x}+\im t \frac{\partial|{\bf \xi}|_\s^2}{\partial \bf x}\right)^a.
$$
Then $b_k({\bf y},{\bf \xi},t)=\sum_{l\leq 2k}b_{k,l}({\bf y},{\bf \xi})t^l.$
\end{lemma}
See \cite{gw} or \cite[Chap. 10]{gs} for more details and a proof of this expansion. The other ingredient is a special case of the lemma of stationary phase.
\begin{lemma}[Lemma of stationary quadratic phase]\label{stationary_phase}
Let $A$ be an $n \times n$ nonsingular self-adjoint matrix, $h\in \bbR^+$, and $f \in C_0^{\infty}(\bbR^n)$.  There is a complete asymptotic expansion
\[
\int_{\bbR^n} f(x) e^{\frac{\sqrt{-1}\langle Ax,x \rangle}{2h}} \mathrm{d}x \sim (2 \pi h)^\frac{n}{2} |\det A|^{-\frac{1}{2}}e^{\frac{i\pi}{4} \operatorname{sgn} A} \left( \text{exp} \left(-\frac{\sqrt{-1}h}{2} b(D)\right) f \right)(0)
\]
where $\operatorname{sgn} A$ is the signature of $A$ and $b(D) = {-}\sum b_{ij} \frac{\partial}{\partial x_i} \frac{\partial}{\partial x_j}$ with $B = A^{-1}$.
\end{lemma}

We are now in a position to prove Theorem \ref{asymp_toric}. 
\begin{proof}
Let ${\bf x}=(u,v)$ be action-angle coordinates for $X$, i.e., $u\in \cp$ is the moment map image of ${\bf x}$ and $v\in \bbR^n/{(2\pi\bbZ)}^n=\bbT^n$. Let ${\bf \xi}=(\hat{u},\hat{v})$ be fiber coordinates on $T_{(u,v)}^*X$.  We can write $\mu_{\frac{\alpha}{\hbar}}(\rho)$ in terms of its Schwartz kernel as
$$
\begin{aligned}
\mu_{\frac{\alpha}{\hbar}}(\rho)&=\int_{\bbT^n}e^{-\frac{\sqrt{-1}\theta \cdot\alpha}{\hbar}}\tr(\action(\theta)_*\rho(\hbar^2\Lap)) \mathrm{d}\theta \\
&=\int_{\bbT^n}e^{-\frac{\sqrt{-1}\theta \cdot \alpha}{\hbar}}\int_{\cp\times \bbT^n}K_{\rho,\hbar}(\action(\theta)(u,v),(u,v))\mathrm{d}u\mathrm{d}v\mathrm{d}\theta\\
&=\int_{\bbT^n}e^{-\frac{\sqrt{-1}\theta \cdot \alpha}{\hbar}}\int_{\cp\times \bbT^n}K_{\rho,\hbar}((u,v+\theta),(u,v))\mathrm{d}u\mathrm{d}v\mathrm{d}\theta.
\end{aligned}
$$
By the Schwartz kernel asymptotic expansion, this expression is $(2\pi\hbar)^{-2n}$ times
$$
\sum_{k\geq 0}\hbar^k\sum_{l\leq 2k}\int_{\bbT^n\times\cp\times \bbT^n\times \bbR^{2n}}b_{k,l}(u,v,\hat{u},\hat{v})\frac{1}{(\im)^l}\frac{d^l\rho}{d\sigma^l}(|{\bf \xi}|_\s^2)e^{\frac{\im(\hat{v}-\alpha)\ip \theta}{\hbar}}\mathrm{d}\hat{u}\mathrm{d}\hat{v}\mathrm{d}u\mathrm{d}v\mathrm{d}\theta
$$
where $b_{k,l}$ are given by Lemma \ref{lemma:schwartz}. Because the metric is torus invariant, the Laplace operator is as well and this implies that the $b_{k,l}$ do not depend on $v$. By changing variables  {so that} $\xi=(\hat{u},\hat{v}+\alpha)$, our expression becomes $(2\pi\hbar)^{-2n}$ times
$$
\sum_{k\geq 0}\hbar^k\sum_{l\leq 2k}\int_{\bbT^n\times\cp\times \bbT^n\times \bbR^{2n}}b_{k,l}(u,\hat{u},\hat{v}+\alpha)\frac{1}{(\im)^l}\frac{d^l\rho}{d\sigma^l}(|{\bf \xi}|_\s^2)e^{\frac{\im\hat{v}\ip \theta}{\hbar}}\mathrm{d}\hat{u}\mathrm{d}\hat{v}\mathrm{d}u\mathrm{d}v\mathrm{d}\theta.
$$

For each $(u,\hat{u},v)$, we are going to apply Lemma \ref{stationary_phase} to the above integral in $(\hat{v},\theta)$ exactly as in \cite[Thm. 5.1]{dgs3}. We can take the matrix $A$ in Lemma \ref{stationary_phase} to be the $2n \times 2n$ matrix given by
$
A=\begin{bmatrix}
0 & {I} \\

{I}&  0
\end{bmatrix}.
$
We have $|\det A |^{-\frac{1}{2}} = 1$, $\sgn A = 0$, and
$
B=\begin{bmatrix}
0 & I \\
I & 0
\end{bmatrix}
$.  Although $\theta$ takes values in $\bbT^n$ and not in $\bbR^n$ as in Lemma \ref{stationary_phase}, the above integral will concentrate on the set where $\theta$ and $\hat{v}$ are zero so that the $\theta$'s outside a square will not contribute mod $O(\hbar^\infty)$. Since the functions $b_{k,l}(u,\hat{u},\hat{v}+\alpha)$ do not depend on $\theta$, we see that applying $b(D)$ to them gives $0$.
Thus we see that  $\mu_{\frac{\alpha}{\hbar}}(\rho)$ is given mod $O(\hbar^\infty)$ by
$$
\mu_{\frac{\alpha}{\hbar}}(\rho)=(2\pi\hbar)^{-n}\sum_{k\geq 0}\hbar^k\sum_{l\leq 2k}\int_{\cp\times  {\bbT^n} \times \bbR^n}b_{k,l}(u,\hat{u},\alpha)\frac{1}{(\im)^l}\frac{d^l\rho}{d\sigma^l}(|{\bf \xi}|_\s^2) {\mathrm{d}\hat{u}}\mathrm{d}u \mathrm{d}v 
$$
where $|{\bf \xi}|_\s^2$ is taken at points where $\xi$ has coordinates $(\hat{u},\alpha)$. Since $b_{k,l}$ does not depend on $v$ and we may take the volume of $\bbT^n=\bbR^n/(2\pi\bbZ)^n$ to be $(2\pi)^n$, we write this as
$$
\hbar^{-n}\sum_{k\geq 0}\hbar^k\sum_{l\leq 2k}\int_{\cp\times \bbR^n}b_{k,l}(u,\hat{u},\alpha)\frac{1}{(\im)^l}\frac{d^l\rho}{d\sigma^l}(|{\bf \xi}|_\s^2)\mathrm{d}u\mathrm{d}\hat{u}. 
$$
At a point $(u,v)$ the norm of the cotangent vector ${\bf \xi}=(\hat{u},\hat{v})$ is
$$
|{\bf \xi}|_\s^2=|(\hat{u},\hat{v})|_\s^2=\hat{u}^t\hess^{-1}(\s)\hat{u}+\hat{v}^t\hess(\s)\hat{v},
$$
and therefore when restricted to the set where $\hat{v}=\alpha$ this gives
\begin{equation}\label{symbol}
|{\bf \xi}|_\s^2=\hat{u}^t\hess^{-1}(\s)\hat{u}+\alpha^t\hess(\s)\alpha.
\end{equation}

Set $b_\cdot^\alpha(u,\hat{u})=b_\cdot(u,\hat{u},\alpha)$. Because the functions $b_k(u,\hat{u},\alpha)$ do not depend on $v$ or $\hat{v}$, the formula for $b_k(u,\hat{u},\alpha)$ gives
\begin{equation}\label{eqn:b_m}
\frac{1}{\im}\frac{\partial b^\alpha_k}{\partial t}=\sum_{\substack{a=(a_1,\ldots, a_{n})\in (\bbZ_+)^{n}\\  |a|\geq 1}}\sum_{j+|a|=k}D^a_{\hat{u}}(|{\bf \xi}|_\s^2)Q_ab_j
\end{equation}
where
$$
Q_a=\frac{1}{a!}\left(\frac{\partial}{\partial u}+\im t \frac{\partial|{\bf \xi}|_\s^2}{\partial u}\right)^a
$$
and $b_k^\alpha(u,\hat{u},t)=\sum_{l\leq 2k}b_{k,l}^\alpha(u,\hat{u})t^l$. From formula (\ref{symbol}) for $|{\bf \xi}|_\s^2$ we see that $D^a_{\hat{u}}(|{\bf \xi}|_\s^2)=D^a_{\hat{u}}(\hat{u}^t\hess^{-1}(\s)\hat{u})$ is zero whenever $|a|>2$ and 
$$
D^a_{\hat{u}}(|{\bf \xi}|_\s^2)=\begin{cases}
\frac{2}{\sqrt{-1}}\sum_{j=1}^n\s^{ij}\hat{u}_j,  \qquad\mbox{if $a$ has a $1$ in position $i$ and $0$ elsewhere;}\\
 -2 \s^{ij},  \qquad \mbox{if $a$ has a $1$ in positions $i$ and $j$ and $0$ elsewhere.}
 \end{cases}
 $$
Note also that if $a$ has a $1$ in position $i$ and $0$ elsewhere then
$$
Q_a=\left(\frac{\partial}{\partial u_i}+\im t \left( \hat{u}^t\frac{\partial\hess^{-1}(\s)}{\partial u_i}\hat{u}+\alpha^t\frac{\partial\hess(\s)}{\partial u_i}\alpha\right)\right),
$$
and if $a$ has a $1$ in positions $i$ and $j$ and $0$ elsewhere then
$$
\begin{aligned}
Q_a=&\frac{1}{2}\left(\frac{\partial}{\partial u_i}+\im t \left( \hat{u}^t\frac{\partial\hess^{-1}(\s)}{\partial u_i}\hat{u}+\alpha^t\frac{\partial\hess(\s)}{\partial u_i}\alpha\right)\right)\cdot\\
&\left(\frac{\partial}{\partial u_j}+\im t \left( \hat{u}^t\frac{\partial\hess^{-1}(\s)}{\partial u_j}\hat{u}+\alpha^t\frac{\partial\hess(\s)}{\partial u_j}\alpha\right)\right)
\end{aligned}
$$
Substituting these expressions into (\ref{eqn:b_m}) gives the expansion in Theorem \ref{asymp_toric}.
\end{proof}

We now make Theorem \ref{asymp_toric} precise in the case of $S^2$.
\begin{thm}\label{two_invariants}
Let $S^2$ be equipped with action-angle coordinates $(x, \theta)$, normalized so that $x \in (-1,1)$.  Consider an $S^1$-invariant metric on $S^2$ given in these coordinates by 
\[
\ddot{\s}dx \otimes dx + \frac{d \theta \otimes d \theta}{\ddot{\s}} .
\]
Let $\alpha\ne 0$ be a real number. Then, for any compactly supported smooth function $\rho$ on $\bbR$, 
\begin{IEEEeqnarray}{c}
\label{inv1}
\int_{[-1,1]\times \mathbb{R}}\rho(\tau) \mathrm{d}x\mathrm{d}\xi
\end{IEEEeqnarray}
and
\begin{IEEEeqnarray}{l}
\label{inv2}
\frac{1}{2} \int_{[-1,1]\times \mathbb{R}} \frac{1}{v} \left[\xi^2 \left(\frac{v^{(2)}}{v^2}-\frac{2(v')^2}{v^3}\right)-\alpha^2 v^{(2)}\right] \rho^{(2)}(\tau) \mathrm{d}x\mathrm{d}\xi \nonumber\\
-\frac{2}{3} \int_{[-1,1]\times \mathbb{R}} \frac{\xi^2}{v} \left[\xi^2 \left(\frac{3(v')^2}{v^4}-\frac{v^{(2)}}{v^3}\right)+\alpha^2\left(\frac{v^{(2)}}{v}-\frac{(v')^2}{v^2}\right)\right] \rho^{(3)}(\tau) \mathrm{d}x\mathrm{d}\xi \nonumber
\\
-\frac{1}{3} \int_{[-1,1]\times \mathbb{R}} \frac{(v')^2}{v} \left(-\frac{\xi^2}{v^2}+\alpha^2\right)^2 \rho^{(3)}(\tau) \mathrm{d}x\mathrm{d}\xi \nonumber
\\
-\frac{1}{2} \int_{[-1,1]\times \mathbb{R}} \frac{\xi^2(v')^2}{v^2} \left(-\frac{\xi^2}{v^2}+\alpha^2\right)^2 \rho^{(4)}(\tau) \mathrm{d}x\mathrm{d}\xi
\end{IEEEeqnarray}
are aespectrally determined, where $v=\ddot{g}$ and $\tau = \frac{\xi^2}{v}+\alpha^2 v$.
\end{thm}
In \cite{dgs3}, the expression in (\ref{inv1}) is shown to be aespectrally determined. Note that the metric defined by $v=\ddot{g}$ is smooth at the poles if and only if $v-\frac{1}{1-x^2}$ is smooth on $[-1,1]$.
\begin{proof}
We will show that the above quantities are the first two nonzero coefficients in the asymptotic expansion of $\mu_{\frac{\alpha}{\hbar}}(\rho)$. The result will follow from the fact that the asymptotic expansion of $\mu_{\frac{\alpha}{\hbar}}$ is aespectrally determined.  By Theorem \ref{asymp_toric}, the spectral measure $\mu_{\frac{\alpha}{\hbar}}$ can be expanded in powers of $\hbar$ as
$$
\hbar^{-1}\sum_{k\geq0} \hbar^k \sum_{l\leq 2k}\int_{[-1,1]\times \bbR}b_{k,l}^\alpha(x,\xi)\frac{1}{(\im)^l}\frac{d^l\rho}{d\sigma^l}\left( \frac{\xi^2}{v}+\alpha^2 v\right)\mathrm{d}x\mathrm{d}\xi
$$
with the functions $b^\alpha_{k,l}$ defined as in the theorem. Since $b^\alpha_0(x,\xi,t)=1$ we see that the leading order term above is simply
$$
\hbar^{-1}\int_{[-1,1]\times \mathbb{R}}\rho\left(\alpha^2 v + \frac{\xi^2}{v}\right) \mathrm{d}x \mathrm{d}\xi.
$$

In our current setting, we note that
\begin{equation}\label{eqn:normS2}
{\hat{u}^t\frac{\partial\hess^{-1}(\s)}{\partial u_i}\hat{u}+\alpha^t\frac{\partial\hess(\s)}{\partial u_i}\alpha = v'\left(\alpha^2 - \frac{\xi^2}{v^2}\right).}
\end{equation}
Thus
$$
\frac{\partial b^\alpha_{1}}{\partial t}=-\frac{2tv'\xi}{v\im} \left(\alpha^2-\frac{\xi^2}{v^2}\right),
$$
so that
$$
b^\alpha_{1,0}=b^\alpha_{1,1}=0, \qquad b^\alpha_{1,2}=-\frac{v'\xi}{v\im} \left(\alpha^2-\frac{\xi^2}{v^2}\right).
$$
We see that the $0$th order term in the expansion of $\mu_{\frac{\alpha}{\hbar}}$ is zero because
$$
\int_{ \bbR}\frac{v'\xi}{v} \left(\alpha^2-\frac{\xi^2}{v^2}\right)\frac{d^2\rho}{d\sigma^2}\left( \frac{\xi^2}{v}+\alpha^2 v\right)\mathrm{d}\xi=0  \text{ for all } x\in (-1,1).
$$

Next we calculate $b_2$. It follows from Theorem \ref{inductionb_m_general} and \eqref{eqn:normS2} that
$$
\begin{aligned}
&\frac{1}{\im}\frac{\partial b^\alpha_{2}}{\partial t}=\frac{2\xi}{\im v}\left(\frac{\partial}{\partial x}+\im t  v'\left(\alpha^2-\frac{\xi^2}{v^2}\right)\right)b^\alpha_{1}\\
&-\frac{1}{v}\left(\frac{\partial}{\partial x}+\im t  v'\left(\alpha^2-\frac{\xi^2}{v^2}\right)\right)^2 {1}
\end{aligned}
$$
This gives 
$$
\begin{aligned}
b^\alpha_2=&\frac{t^2}{2v}\left( v'\left(\alpha^2-\frac{\xi^2}{v^2}\right) \right)'\\
&+\frac{\im t^3}{3}\left( \frac{(v')^2}{v}\left(\alpha^2-\frac{\xi^2}{v^2}\right)^2+\frac{2\xi^2}{v}\left( \frac{v'}{v}\left(\alpha^2-\frac{\xi^2}{v^2}\right) \right)'\right)\\
&-\frac{t^4(v')^2\xi^2}{2v^2}\left(\alpha^2-\frac{\xi^2}{v^2}\right)^2 \\
\end{aligned}
$$
Therefore the coefficient of $\hbar$ in the expansion of $\mu_{\frac{\alpha}{\hbar}}$ is 
\begin{IEEEeqnarray*}{l}
\int_{[-1,1]\times \bbR}\left[\frac{1}{2v}\left( \xi^2 \left( \frac{v^{(2)}}{v^2}-\frac{2(v')^{2}}{v^3}\right)-\alpha^2 v^{(2)}\right) \frac{d^2\rho}{d\sigma^2}\left( \frac{\xi^2}{v}+\alpha^2 v\right)\right.\\
-\frac{(v')^{2}}{3v}\left(\alpha^2-\frac{\xi^2}{v^2}\right)^2\frac{d^3\rho}{d\sigma^3}\left( \frac{\xi^2}{v}+\alpha^2 v\right)\\
-\frac{2\xi^2}{3v}\left( \xi^2 \left( -\frac{v^{(2)}}{v^3}+\frac{3(v')^{2}}{v^4}\right)+\alpha^2  \left( \frac{v^{(2)}}{v}-\frac{(v')^{2}}{v^2}\right)\right) \frac{d^3\rho}{d\sigma^3}\left( \frac{\xi^2}{v}+\alpha^2 v\right)\\
\left.-\frac{(v')^2\xi^2}{2v^2}\left(\alpha^2-\frac{\xi^2}{v^2}\right)^2 \frac{d^4\rho}{d\sigma^4}\left( \frac{\xi^2}{v}+\alpha^2 v\right)\right]\mathrm{d}x\mathrm{d}\xi
\end{IEEEeqnarray*}
and the result follows.
\end{proof}


\section{Inverse spectral results for $S^1${-}invariant metrics on $S^2$}\label{proof}

In this section we prove {Theorem \ref{main}. In \cite[Thm. 6.15]{dgs3}, the first invariant in Theorem \ref{two_invariants} was used to show that when $v$ is even and convex,} the {corresponding $S^1$-invariant metric on $S^2$} is determined by the asymptotic equivariant spectrum. Colin de Verdi\`ere {\cite{cdv}} has shown {that} the spectrum of a Schr\"{o}dinger operator on $\mathbb{R}^2$ with a {\em single well} potential essentially determines the potential. In \cite{gw} Guillemin and Wang use higher order semi-classical spectral invariants to generalize this result to {\em double wells}. {Our approach} is analogous to the Guillemin-Wang generalization of Colin de Verdi\`ere's result. Using the {new spectral invariant in Theorem \ref{two_invariants}, we generalize Theorem 6.15 in \cite{dgs3}} to show that when $v$ is a single well, the asymptotic equivariant spectrum determines $v$ and hence the metric. {Functions that are even and convex are clearly} very special cases of single wells.

Suppose now that $v$ is a single well, i.e., suppose it has a unique nondegenerate minimum at $x=0$, and that $v$ is increasing for $x$ positive and decreasing for $x$ negative. Let $c\neq 0$ be the {minimum} value of $v$ at the point $x=0$. Let $\alpha\ne 0$ be a fixed real number. We will show how to use invariants (\ref{inv1}) and (\ref{inv2}) to recover the function $v(x)$ on the interval $|x|<1$.   {Note that these invariants involve integrals of the function $\rho(\tau) = \rho\left(\frac{\xi^2}{v(x)}+\alpha^2 v(x)\right)$ and its derivatives.  To evaluate such an integral using} Fubini's theorem over the region 
$$
\{ (x,\xi): x\in [-1,1], \frac{\xi^2}{v(x)}+\alpha^2v(x)<\lambda  \},
$$
one needs to break up the integral into several integrals corresponding to regions where the condition $\frac{\xi^2}{v(x)}+\alpha^2v(x)=\lambda$ defines $\xi$ as a function of $x$. The single well condition ensures that we only need to  {consider two regions; this is the only way in which} we use this condition. For $0<\lambda <1$, we denote by $A_1^{\lambda}$ the region in the first quadrant bounded by the curve $\frac{\xi^2}{v(x)}+\alpha^2 v(x)=\lambda$ and by $A_2^{\lambda}$ the region in the second quadrant bounded by the same curve.  {Thus} $A_1^{\lambda}$ can be described as
\[
A_1^{\lambda} = \{ (x,\xi): x\in [0,1], \quad v(x)< \frac{\lambda}{\alpha^2},\quad 0<\xi<\sqrt{\lambda v(x) - \alpha^2 v(x)^2}  \}
\]
and $A_2^{\lambda}$ can be described as
\[
A_2^{\lambda} = \{ (x,\xi): x\in [-1,0], \quad v(x)< \frac{\lambda}{\alpha^2},\quad 0<\xi<\sqrt{\lambda v(x) - \alpha^2 v(x)^2}   \}.
\]
Note that both sets can be empty for a given $\lambda$. Since $v$ is a single well with minimum at $0$,  {the regions} $A_1^{\lambda}$ and $A_2^{\lambda}$ are, indeed, bounded.   {If we could choose $\rho_{\lambda}$ to be the characteristic function of $[0,\lambda]$, we would see from invariant (\ref{inv1}) that the {sum}
\begin{IEEEeqnarray}{c}
\label{exp:firstint}
\int_{A_1^\lambda}{\mathrm{d}x \mathrm{d}\xi} + \int_{A_2^\lambda} \mathrm{d}x \mathrm{d}\xi
\end{IEEEeqnarray}
would be aespectrally determined. However, such a $\rho_\lambda$ is not smooth, so to make this precise we must consider $\rho_\lambda$ as the limit of an appropriate sequence of functions that equal the characteristic function of $[0,\lambda]$ on larger and larger subsets of $[0,\lambda]$.}

Let $x=f_1(s)$ be the inverse function of $s=v({x})$, for $x \in (0,1)$. Then
\begin{IEEEeqnarray}{rCl}
\int_{A_1^\lambda} \mathrm{d}x \mathrm{d}\xi & = & \int_0^{f_1\left(\frac{\lambda}{\alpha^2}\right)} \int_0^{\sqrt{\lambda v-\alpha^2 v^2}} \mathrm{d}\xi \mathrm{d}x \nonumber \\
& = & \int_0^{f_1\left(\frac{\lambda}{\alpha^2}\right)} \sqrt{\lambda v(x)-\alpha^2 v(x)^2} \mathrm{d}x \nonumber \\
& = & \int_c^{\frac{\lambda}{\alpha^2}} \sqrt{\lambda s - \alpha^2 s^2} \frac{df_1}{ds} \mathrm{d}s . \label{exp1}
\end{IEEEeqnarray}
{Analogously, we define $x=f_2(s)$ to be the inverse function of $s=v(-x)$, for $x \in (0,1)$.  Then the same calculations give}
\begin{IEEEeqnarray}{rCl}
\label{exp2}
\int_{A_2^\lambda} \mathrm{d}x \mathrm{d}\xi & = & \int_c^{\frac{\lambda}{\alpha^2}} \sqrt{\lambda s - \alpha^2 s^2} \frac{df_2}{ds} \mathrm{d}s {.}
\end{IEEEeqnarray}
This implies in particular that $c$ is aespectrally determined{:} if $\lambda/\alpha^2<c$, then integral (\ref{exp:firstint}) is zero{,} whereas it is non-zero if $\lambda/\alpha^2>c$.

Now set $S=s-c$ and $\beta=\frac{\lambda}{\alpha^2}-c$. {Substituting the expressions in \eqref{exp1} and \eqref{exp2} into \eqref{exp:firstint}, we see that the} following quantity is aespectrally determined for all $\beta$:
\begin{equation}\label{abel_invt1}
\int_0^{\beta} \sqrt{\beta-S}   {\sqrt{S+c}\left(\left(\frac{df_1}{dS}+\frac{df_2}{dS} \right)\bigg|_{S+c}\right)} \mathrm{d}S {.}
\end{equation}
{The function in \eqref{abel_invt1} can be viewed as the Abel transform of another function, as we now explain.}  Recalling that the fractional integra{tion} operation of Abel is defined as
\begin{IEEEeqnarray*}{c}
J^a g(s)=\frac{1}{\Gamma(a)} \int_0^s (s-\nu)^{a-1} g(\nu)\mathrm{d}\nu,
\end{IEEEeqnarray*}
for $a>0$, we observe that (\ref{abel_invt1}) corresponds to
\begin{IEEEeqnarray*}{c}
\Gamma \left(\frac{3}{2}\right) J^{\frac{3}{2}}\left(\sqrt{S+c}   {\left(\left(\frac{df_1}{dS}+\frac{df_2}{dS} \right)\bigg|_{S+c}\right)}\right)(\beta) {.}
\end{IEEEeqnarray*}
As the Abel transform of a function determines the function, we may recover the quantity $\sqrt{S+c}   {\left(\left(\frac{df_1}{dS}+\frac{df_2}{dS} \right)\bigg|_{S+c}\right)}$ and hence $\frac{df_1}{dS}+\frac{df_2}{dS}$ as a function of $s$. For more on the Abel transform and its invertibility see {\cite[\S 10.6]{gs}.}

Next we integrate the first and last terms of invariant (\ref{inv2}). Integration by parts with respect to $\xi$ gives
\begin{eqnarray*}
&&\frac{1}{2} \int_{[-1 ,1] \times \mathbb{R}} \frac{1}{v} \left[\xi^2 \left(\frac{v^{(2)}}{v^2}-\frac{2(v')^2}{v^3}\right)-\alpha^2 v^{(2)}\right] \rho^{(2)}(\tau) \mathrm{d}x\mathrm{d}\xi \\
&&=\frac{1}{2} \int_{[-1 ,1] \times \mathbb{R}} \frac{1}{v}\frac{\partial}{\partial \xi} \left[\frac{\xi^3}{3} \left(\frac{v^{(2)}}{v^2}-\frac{2(v')^2}{v^3}\right)-\xi\alpha^2 v^{(2)}\right] \rho^{(2)}(\tau) \mathrm{d}x\mathrm{d}\xi \\
&&=-\frac{1}{2} \int_{[-1 ,1] \times \mathbb{R}} \left[\frac{2\xi^4}{3v^2} \left(\frac{v^{(2)}}{v^2}-\frac{2(v')^2}{v^3}\right)-\frac{2\xi^2\alpha^2 v^{(2)}}{v^2}\right] \rho^{(3)}(\tau) \mathrm{d}x\mathrm{d}\xi,
\end{eqnarray*}
and
\begin{eqnarray*}
&&-\frac{1}{2} \int_{[-1 ,1] \times \mathbb{R}} \frac{\xi^2(v')^2}{v^2} \left(-\frac{\xi^2}{v^2}+\alpha^2\right)^2 \rho^{(4)}(\tau) \mathrm{d}x\mathrm{d}\xi \\
&&=-\frac{1}{2} \int_{[-1 ,1] \times \mathbb{R}} \frac{\xi(v')^2}{2{v}} \left(-\frac{\xi^2}{v^2}+\alpha^2\right)^2\frac{\partial( \rho^{(3)}(\tau))}{\partial \xi} \mathrm{d}x\mathrm{d}\xi \\
&&=\frac{1}{2} \int_{[-1 ,1] \times \mathbb{R}} \left[\frac{5\xi^4 (v')^2}{2 v^5}-\frac{3\xi^2 \alpha^2 (v')^2}{v^3} + \frac{\alpha^4 (v')^2}{2v}\right] \rho^{(3)}(\tau) \mathrm{d}x\mathrm{d}\xi.
\end{eqnarray*}
Note that we do not pick up boundary terms because $\rho$ is compactly supported. {Combining these new expressions for the first and last terms of \eqref{inv2} with the middle terms of \eqref{inv2}, we} conclude that
\begin{IEEEeqnarray}{l}
\label{terceira_deriv}
\int_{[-1,1]\times \bbR}
 \left[ -\frac{1}{2} \left(\frac{2\xi^4}{3v^2}\left(\frac{v^{(2)}}{v^2}-\frac{2(v')^2}{v^3} \right)-\frac{2\alpha^2 v^{(2)}\xi^2}{v^2} \right) \nonumber \right. \\
-\frac{2\xi^2}{3v}\left(\xi^2\left(\frac{3(v')^2}{v^4}-\frac{v^{(2)}}{v^3}\right) + \alpha^2\left(\frac{v^{(2)}}{v}-\frac{(v')^2}{v^2}\right) \right) \nonumber \\
-\frac{(v')^2}{3v} \left(-\frac{\xi^2}{v^2}+\alpha^2\right)^2 \nonumber \\
+\frac{(v')^2}{2v} \left. \left(\frac{5\xi^4}{2v^4}-\frac{3\xi^2\alpha^2}{v^2}+\frac{\alpha^4}{2}\right) \right] \rho^{(3)}(\tau)\mathrm{d}x\mathrm{d}\xi
\end{IEEEeqnarray}
is aespectrally determined for all $\alpha$ and all compactly supported $\rho$.  

By taking {the limit of an appropriate sequence of functions that equal $\rho$ on larger and larger sets, we see that we can make $\rho_\Lambda(\tau)=e^{-\Lambda \tau}$ in (\ref{terceira_deriv}). It} follows that
\begin{IEEEeqnarray*}{l}
\int_{[-1,1]\times \bbR}
 \left[ -\frac{1}{2} \left(\frac{2\xi^4}{3v^2}\left(\frac{v^{(2)}}{v^2}-\frac{2(v')^2}{v^3} \right)-\frac{2\alpha^2 v^{(2)}\xi^2}{v^2} \right) \nonumber \right. \\
-\frac{2\xi^2}{3v}\left(\xi^2\left(\frac{3(v')^2}{v^4}-\frac{v^{(2)}}{v^3}\right) + \alpha^2\left(\frac{v^{(2)}}{v}-\frac{(v')^2}{v^2}\right) \right) \nonumber \\
-\frac{(v')^2}{3v} \left(-\frac{\xi^2}{v^2}+\alpha^2\right)^2 \nonumber \\
+\frac{(v')^2}{2v} \left. \left(\frac{5\xi^4}{2v^4}-\frac{3\xi^2\alpha^2}{v^2}+\frac{\alpha^4}{2}\right) \right] \Lambda^3 e^{-\Lambda \tau}\mathrm{d}x\mathrm{d}\xi  \\
\end{IEEEeqnarray*}
is aespectrally determined for all $\alpha$ and $\Lambda${; thus}
\begin{IEEEeqnarray*}{l}
\int_{[-1,1]\times \bbR}
 \left[ -\frac{1}{2} \left(\frac{2\xi^4}{3v^2}\left(\frac{v^{(2)}}{v^2}-\frac{2(v')^2}{v^3} \right)-\frac{2\alpha^2 v^{(2)}\xi^2}{v^2} \right) \nonumber \right. \\
-\frac{2\xi^2}{3v}\left(\xi^2\left(\frac{3(v')^2}{v^4}-\frac{v^{(2)}}{v^3}\right) + \alpha^2\left(\frac{v^{(2)}}{v}-\frac{(v')^2}{v^2}\right) \right) \nonumber \\
-\frac{(v')^2}{3v} \left(-\frac{\xi^2}{v^2}+\alpha^2\right)^2 \nonumber \\
+\frac{(v')^2}{2v} \left. \left(\frac{5\xi^4}{2v^4}-\frac{3\xi^2\alpha^2}{v^2}+\frac{\alpha^4}{2}\right) \right] e^{-\Lambda \tau}\mathrm{d}x\mathrm{d}\xi  \\
\end{IEEEeqnarray*}
is aespectrally determined for all $\alpha$ and $\Lambda$. {We may} approximate the characteristic function of $[0,\lambda]$ by a linear combination of functions {of the form} $e^{-\Lambda \tau}$ where $\Lambda$ may be complex{, implying that 
\begin{IEEEeqnarray}{l}
\label{inv2a}
\int_{A_1^\lambda+A_2^\lambda} \left[-\frac{1}{2} \left(\frac{2\xi^4}{3v^2}\left(\frac{v^{(2)}}{v^2}-\frac{2(v')^2}{v^3} \right)-\frac{2\alpha^2 v^{(2)}\xi^2}{v^2} \right) \right.\nonumber \\
-\frac{2\xi^2}{3v}\left(\xi^2\left(\frac{3(v')^2}{v^4}-\frac{v^{(2)}}{v^3}\right) + \alpha^2\left(\frac{v^{(2)}}{v}-\frac{(v')^2}{v^2}\right) \right) \nonumber \\
-\frac{(v')^2}{3v} \left(-\frac{\xi^2}{v^2}+\alpha^2\right)^2 \nonumber \\
+\left. \frac{(v')^2}{2v} \left(\frac{5\xi^4}{2v^4}-\frac{3\xi^2\alpha^2}{v^2}+\frac{\alpha^4}{2}\right)\right] \mathrm{d}x\mathrm{d}\xi
\end{IEEEeqnarray}
is aespectrally determined for all $\alpha$ and $\lambda$.}

We will treat the integral over $A_1^\lambda$ first. By integrating the part of expression (\ref{inv2a}) concerning the region $A_1^\lambda$ with respect to $\xi$, we obtain
\begin{IEEEeqnarray*}{l}
- \int_0^{f_1\left(\frac{\lambda}{\alpha^2}\right)} \frac{1}{2}\left[\frac{2\xi^5}{15v^2}\left(\frac{v^{(2)}}{v^2}-\frac{2(v')^2}{v^3}\right)-\frac{2\alpha^2 v^{(2)} \xi^3}{3v^2}\right]_0^{\sqrt{\lambda v-\alpha^2 v^2}} \mathrm{d}x\\
- \int_0^{f_1\left(\frac{\lambda}{\alpha^2}\right)} \frac{2}{3}\left[\frac{\xi^5}{5v}\left(\frac{3(v')^2}{v^4}-\frac{v^{(2)}}{v^3}\right)+\frac{\xi^3 \alpha^2}{3v}\left(\frac{v^{(2)}}{v}-\frac{(v')^2}{v^2}\right) \right]_0^{\sqrt{\lambda v-\alpha^2 v^2}} \mathrm{d}x \\
- \int_0^{f_1\left(\frac{\lambda}{\alpha^2}\right)} \frac{1}{3}\left[\frac{\xi^5 (v')^2}{5v^5}-\frac{2\xi^3 \alpha^2 (v')^2}{3v^3}+\frac{(v')^2 \alpha^4 \xi}{v} \right]_0^{\sqrt{\lambda v-\alpha^2 v^2}} \mathrm{d}x \\
+ \int_0^{f_1\left(\frac{\lambda}{\alpha^2}\right)} \frac{1}{2}\left[\frac{\xi^5(v')^2}{2v^5}-\frac{\xi^3 \alpha^2 (v')^2}{v^3}+\frac{\alpha^4 (v')^2 \xi}{2v}\right]_0^{\sqrt{\lambda v-\alpha^2 v^2}} \mathrm{d}x.
\end{IEEEeqnarray*}
{A}fter rearranging terms, we get
\begin{IEEEeqnarray*}{c}
\int_0^{f_1\left(\frac{\lambda}{\alpha^2}\right)} \left[-\frac{\alpha^4}{12} \frac{(v')^2 \xi}{v}+\frac{\alpha^2}{9}\frac{\xi^3 v^{(2)}}{v^2}-\frac{\alpha^2}{18} \frac{\xi^3 (v')^2}{v^3}+ \frac{1}{15}\frac{\xi^5 v^{(2)}}{v^4}-\frac{1}{12}\frac{\xi^5 (v')^2}{v^5}\right]_0^{\sqrt{\lambda v-\alpha^2 v^2}} \mathrm{d}x {,}
\end{IEEEeqnarray*}
{or equivalently,}
\begin{IEEEeqnarray}{c}
\label{exp:changevar}
\int_0^{f_1\left(\frac{\lambda}{\alpha^2}\right)} \sqrt{\lambda v-\alpha^2 v^2}\left(-\frac{\alpha^4}{9}\frac{(v')^2}{v}-\frac{\lambda \alpha^2}{45}\frac{v^{(2)}}{v} +\frac{\lambda \alpha^2}{9}\frac{(v')^2}{v^2}\right. \nonumber \\
\left. + \frac{\lambda^2}{15}\frac{v^{(2)}}{v^2}-\frac{2\alpha^4 v^{(2)}}{45} -\frac{\lambda^2}{12}\frac{(v')^2}{v^3}\right) \mathrm{d}x {.}
\end{IEEEeqnarray}

Recalling that $x=f_1(s)$ is the inverse function of $s=v(x)$ for $x \in (0,{1}) $, we have
\begin{IEEEeqnarray*}{l}
v'(f_1(s))=\frac{1}{f'_1(s)},\\
v''(f_1(s)) = -\frac{f''_1(s)}{(f'_1(s))^3}.
\end{IEEEeqnarray*}
Expression (\ref{exp:changevar}) can now be rewritten as
\begin{eqnarray}
\int_c^{\frac{\lambda}{\alpha^2}} \sqrt{\lambda-\alpha^2 s}\left[-\frac{\alpha^4}{9}\frac{1}{\sqrt{s}}\frac{1}{\frac{df_1}{ds}}+\frac{\lambda \alpha^2}{45}\frac{1}{\sqrt{s}}\frac{\frac{d^2 f_1}{ds^2}}{\left(\frac{df_1}{ds}\right)^2}+\frac{\lambda \alpha^2}{9} \frac{1}{s\sqrt{s}}\frac{1}{\frac{d f_1}{ds}}  \right.\nonumber \\
 \left.-\frac{\lambda^2}{15}\frac{1}{s\sqrt{s}}\frac{\frac{d^2 f_1}{ds^2}}{\left(\frac{df_1}{ds}\right)^2}+\frac{2\alpha^{4}}{45} \sqrt{s} \frac{\frac{d^2 f_1}{ds^2}}{\left(\frac{df_1}{ds}\right)^2}-\frac{\lambda^2}{12}\frac{1}{s^2 \sqrt{s}}\frac{1}{\frac{df_1}{ds}}\right] \mathrm{d}s  \nonumber
\\
=\int_c^{\frac{\lambda}{\alpha^2}}\sqrt{\lambda-\alpha^2 s} \left[\left(-\frac{\alpha^4}{9}\frac{1}{\sqrt{s}}+\frac{\lambda \alpha^2}{9}\frac{1}{s\sqrt{s}}-\frac{\lambda^2}{12}\frac{1}{s^2\sqrt{s}}\right)\frac{1}{\frac{df_1}{ds}} \right.\nonumber \\
 \left.+ \left(\frac{\lambda \alpha^2}{45}\frac{1}{\sqrt{s}}-\frac{\lambda^2}{15}\frac{1}{s\sqrt{s}}+\frac{2 \alpha^{4}}{45} \sqrt{s}\right) \frac{\frac{d^2 f_1}{ds^2}}{\left(\frac{df_1}{ds}\right)^2}\right] \mathrm{d}s . \label{exp:beforechange}
\end{eqnarray}

Next we make a change of variable, setting $\beta=\frac{\lambda}{\alpha^2}-c$ and $S=s-c$, so that (\ref{exp:beforechange}) can be rewritten as
\begin{IEEEeqnarray*}{c}
 \alpha^5 \int_0^{\beta} \sqrt{\beta-S}\left[\left(-\frac{1}{9\sqrt{S+c}}+\frac{\beta+c}{9(S+c)\sqrt{S+c}}-\frac{(\beta+c)^2}{12(S+c)^2\sqrt{S+c}}\right)
    { \left(\frac{1}{\frac{df_1}{dS}}\right)\bigg|_{S+c} }   
  \right. \nonumber \\
 \left. + \left(\frac{\beta+c}{45\sqrt{S+c}}-\frac{(\beta+c)^2}{15(S+c)\sqrt{S+c}}+\frac{2\sqrt{S+c}}{45}\right)   {\left( \frac{\frac{d^2 f_1}{d S^2}}{\left(\frac{df_1}{dS}\right)^2}\right)\bigg|_{S+c}}\right] \mathrm{d}S .
\end{IEEEeqnarray*}
We {repeat} this calculation for $A_2^\lambda$ and we conclude that the second order spectral invariant \eqref{inv2a} is given by
\begin{IEEEeqnarray*}{c}
 \alpha^5 \int_0^{\beta} \sqrt{\beta-S}\left[\left(-\frac{1}{9\sqrt{S+c}}+\frac{\beta+c}{9(S+c)\sqrt{S+c}}-\frac{(\beta+c)^2}{12(S+c)^2\sqrt{S+c}}\right)\right. \\
 \left.   {\left( \frac{1}{\frac{df_1}{dS}}+\frac{1}{\frac{df_2}{dS}}\right)\bigg|_{S+c} }\right. \\
\left. + \left(\frac{\beta+c}{45\sqrt{S+c}}-\frac{(\beta+c)^2}{15(S+c)\sqrt{S+c}}+\frac{2\sqrt{S+c}}{45}\right){ \left(\frac{\frac{d^2 f_1}{dS^2}}{\left(\frac{df_1}{dS}\right)^2}+\frac{\frac{d^2 f_2}{dS^2}}{\left(\frac{df_2}{dS}\right)^2} \right)\bigg|_{S+c}}\right] \mathrm{d}S . \\
\end{IEEEeqnarray*}
Using Abel's fractional integra{tion}, we observe that {the preceding} expression corresponds to
\begin{IEEEeqnarray*}{c}
\Gamma \left(\frac{3}{2}\right) J^{\frac{3}{2}}\left[ \left(-\frac{1}{9\sqrt{S+c}}+\frac{\beta+c}{9(S+c)\sqrt{S+c}}-\frac{(\beta+c)^2}{12(S+c)^2\sqrt{S+c}}\right)\right.\\
\left. \left( \frac{1}{\frac{df_1}{dS}}+\frac{1}{\frac{df_2}{dS}} \right){\bigg|_{S+c}} \right. \\
\left. +\left(\frac{\beta+c}{45\sqrt{S+c}}-\frac{(\beta+c)^2}{15(S+c)\sqrt{S+c}}+{\frac{2\sqrt{S+c}}{45}}\right)\left(\frac{\frac{d^2 f_1}{dS^2}}{\left(\frac{df_1}{dS}\right)^2}+\frac{\frac{d^2 f_2}{dS^2}}{\left(\frac{df_2}{dS}\right)^2} \right){\bigg|_{S+c}}\right] (\beta).\\
\end{IEEEeqnarray*}
As a function {is} uniquely determined by its Abel transform, we recover the quantity
\begin{IEEEeqnarray*}{c}
\left[ \left(-\frac{1}{9\sqrt{S+c}}+\frac{\beta+c}{9(S+c)\sqrt{S+c}}-\frac{(\beta+c)^2}{12(S+c)^2\sqrt{S+c}}\right)\left( \frac{1}{\frac{df_1}{dS}}+\frac{1}{\frac{df_2}{dS}} \right){\bigg|_{S+c}} \right. \\
\left. + \left(\frac{\beta+c}{45\sqrt{S+c}}-\frac{(\beta+c)^2}{15(S+c)\sqrt{S+c}}+{\frac{2\sqrt{S+c}}{45}}\right)\left(\frac{\frac{d^2 f_1}{dS^2}}{\left(\frac{df_1}{dS}\right)^2}+\frac{\frac{d^2 f_2}{dS^2}}{\left(\frac{df_2}{dS}\right)^2} \right){\bigg|_{S+c}}\right].\\
\end{IEEEeqnarray*}
We rewrite {this quantity} as 
\begin{equation}\label{ode}
A_{\beta,c}(S)\left( \frac{1}{\frac{df_1}{dS}}+\frac{1}{\frac{df_2}{dS}} \right) +B_{\beta,c}(S)\frac{d}{dS}\left( \frac{1}{\frac{df_1}{dS}}+\frac{1}{\frac{df_2}{dS}} \right), \\
\end{equation}
where we set
$$
A_{\beta,c}(S)=-\frac{1}{9\sqrt{S+c}}+\frac{\beta+c}{9(S+c)\sqrt{S+c}}-\frac{(\beta+c)^2}{12(S+c)^2\sqrt{S+c}}
$$
and
$$
B_{\beta,c}(S)=-\frac{\beta+c}{45\sqrt{S+c}}+\frac{(\beta+c)^2}{15(S+c)\sqrt{S+c}}-{\frac{2\sqrt{S+c}}{45}}.
$$

Our goal is to show that the function \eqref{ode} completely determines $\left(\frac{1}{\frac{df_1}{dS}}+\frac{1}{\frac{df_2}{dS}}\right)\bigg|_{S+c}$ as this will imply that  $\frac{1}{\frac{df_1}{dS}}+\frac{1}{\frac{df_2}{dS}}$ itself is determined. {S}ince $v$ has a minimum at $0$, {we have $v'(0)=0$ and} the derivatives of $f_1$ and $f_2$ tend to $+\infty$ as $s$ tends to $c$, i.e., as $S$ tends to $0$. This implies that
$$
{\left( \frac{1}{\frac{df_1}{dS}}+\frac{1}{\frac{df_2}{dS}} \right)\bigg|_{S+c}(0)=\left( \frac{1}{\frac{df_1}{dS}}+\frac{1}{\frac{df_2}{dS}} \right)\bigg|_{s=c}=0. }
$$
We want to show that the initial value ODE 
$$
\begin{cases}
A_{\beta,c}(S)F(S)+B_{\beta,c}(S)F'(S)=K(S), \quad S\in [0,\beta] \\
F(0)=0
\end{cases}
$$
has at most one solution, for any given function $K$. Let {$M$ and $N$} be two solutions{.  Then $M-N=f$} satisfies
$$
\begin{cases}
A_{\beta,c}(S)f(S)+B_{\beta,c}(S)f'(S)=0,\quad S\in [0,\beta]\\
f(0)=0.
\end{cases}
$$
We show that $f$ vanishes identically. {Explicit integration of this first-order linear equation gives}
$$
f(S)=Ce^{-\int_0^S\frac{A_{\beta,c}(\tau)}{B_{\beta,c}(\tau)}\mathrm{d}\tau},
$$
for some constant $C$. {One can} check that 
$$
A_{\beta,c}(S)=-\frac{(S+c)^2-(S+c)(\beta+c)+\frac{3}{4}(\beta+c)^2}{9(S+c)^2\sqrt{S+c}}
$$
and
$$
B_{\beta,c}(S)=\frac{({\beta-S})(2S+3\beta+5c)}{45(S+c)\sqrt{S+c}}
$$
so that 
$$
\frac{A}{B}=-\frac{5\left((S+c)^2-(S+c)(\beta+c)+\frac{3}{4}(\beta+c)^2\right)}{(S+c)({\beta-S})(2S+3\beta+5c)},
$$
which we can rewrite as 
$$
\frac{A}{B}=\frac{q_1}{S+c}+\frac{q_2}{{\beta-S}}+\frac{q_3}{2S+3\beta+5c}
$$
for some constants $q_1,q_2, q_3${. Thus} 
$$
\int_0^S\frac{A_{\beta,c}(\tau)}{B_{\beta,c}(\tau)}\mathrm{d}\tau={q_1}\log({S+c}){-}{q_2}\log({{\beta-S}})+{{\frac{q_3}{2}}}\log({2S+3\beta+5c}),
$$
{and hence} 
$$
f(S)=C(S+c)^{{-}q_1}({\beta-S})^{q_2}(2S+3\beta+5c)^{{-\frac{q_3}{2}}}.
$$
Because $f(0)=0$, this forces $C=0$ and therefore $f$ is identically zero as claimed. 

We conclude that
$$
\frac{1}{\frac{d f_1}{d S}}+ \frac{1}{\frac{d f_2}{d S}}
$$
{is} aespectrally determined{.  Recall that we used a straightforward Abel transform argument to show that
$$\frac{df_1}{dS}+\frac{df_2}{dS}$$
is aespectrally determined.  But if we know $\frac{1}{p} + \frac{1}{q}= \frac{p+q}{pq}$ and we know $p+q$, then we clearly know $pq$.  This in turn implies that we know $|p-q|$, since $(p-q)^2 = (p+q)^2-4pq$.  Thus we know $p$ and $q$, up to order.  Hence the functions $\frac{df_1}{dS}$ and $\frac{df_2}{dS}$ are aespectrally determined, and we recover the functions $f_1$ and $f_2$   {up to order, i.e., up to knowing which function is the inverse of $v$ in which quadrant}.  {Therefore the function $v(x)$ on the interval $|x|<  {1}$ is aespectrally determined   {up to reversing the $x$-coordinate on $[-1,1]$.  Note that this reversal gives rise to an isometric metric on $S^2$, so that} we have proved Theorem \ref{main}.


\begin{thebibliography}{99}

\bibitem[A]{a}  M. Abreu, {\em K\"ahler geometry of toric manifolds in symplectic coordinates}, Symplectic and contact topology: interactions and perspectives (Toronto, ON/Montreal, QC, 2001), Fields Inst. Commun., {\bf 35}, (2003), 1--24.

\bibitem[AF]{af} M. Abreu and P. Freitas, {\em On the invariant spectrum of $S^1$-invariant metrics on $S^2$}, Proc. London Math. Soc. (3) {\bf 84} (2002), no. 1, 213--230.

\bibitem[BH]{bh} J. Br{\"u}ning and E. Heintze, {\em Spektrale {S}tarrheit gewisser {D}rehfl\"achen}, Math. Ann. {\bf 269} (1984), no. 1, 95--101.

\bibitem[CdV]{cdv} Y. Colin de Verdi{\`e}re, {\em  A semi-classical inverse problem II: reconstruction of the potential}, Geometric aspects of analysis and mechanics, Progr. Math., {\bf 292} (2011), 97--119.

\bibitem[DGS1]{dgs1} E. Dryden, V. Guillemin, and R. Sena-Dias, {\em Hearing Delzant polytopes from the equivariant spectrum}, Trans. Amer. Math. Soc. {\bf 364} (2012), no. 2, 887--910.

\bibitem[DGS2]{dgs2} E. Dryden, V. Guillemin, and R. Sena-Dias, {\em Equivariant inverse spectral theory and toric orbifolds}, Adv. Math. {\bf 231} (2012), no. 3-4, 1271--1290.

\bibitem[DGS3]{dgs3}  E. Dryden, V. Guillemin, and R. Sena-Dias, {\em Semi-classical weights and equivariant spectral theory}, arXiv:1401.8285.

\bibitem[G]{g} V. Guillemin, {\em Kaehler structures on toric varieties}, J. Differential Geom. {\bf 40} (1994), no. 2, 285--309.

\bibitem[GS]{gs} V.~Guillemin and S. Sternberg,  {\em Semi-classical analysis},  International Press, Boston, MA, 2013.

\bibitem[GW]{gw} V. Guillemin and Z. Wang,  {\em Semiclassical spectral invariants for Schr\"{o}dinger operators}, J. Differential Geom. {\bf 91} (2012), no. 1, 103--128.

\bibitem[KW]{kw} J. Kazdan and F. Warner, {\em Surfaces of revolution with monotonic increasing curvature and an application to the equation $u=1-Ke^{2u}$ on $S^2$}, Proc. Amer. Math. Soc. {\bf 32} (1972),139--141.

\bibitem[M]{m} J. Milnor, {\em Eigenvalues of the Laplace operator on certain manifolds}, Proc. Nat. Acad. Sci. U.S.A. {\bf 51} (1964), 542.

\bibitem[Z]{z} S. Zelditch {\em The inverse spectral problem for surfaces of revolution}, J. Differential Geom. {\bf 49} (1998), 2007--264.


\end{thebibliography}
\end{document}